\documentclass[11pt,notitlepage,oneside,openany]{amsart}
\usepackage{amssymb}
\usepackage{amsthm}
\usepackage{mathpple}

\newcommand{\fract}[2]{\textstyle\frac{{#1}}{{#2}}}

\theoremstyle{plain}
\newtheorem*{thm}{Theorem}
\newtheorem{lem}{Lemma}

\begin{document}
\begin{center}
{\bf\Large On a sequence related to the Josephus problem }\\[5pt]
\textsc{Ralf Stephan\footnote{The author can be reached at
\texttt{mailto:ralf@ark.in-berlin.de}.}}
\end{center}
{\Small\baselineskip6pt
In this short note, we show that an integer sequence defined on the minimum of differences
between divisor complements of its partial products is connected with
the Josephus problem (q=3).
}
%\end{abstract}
%\maketitle
\bigskip

We prove the following theorem and, finally, state the relatedness of two constants.
\begin{thm}
Let $a_n$ and $b_n$ be recursively defined as
\begin{align}
a_0 & = 4,\;
a_{n} = \min\left(|d_j-{p_n}/{d_j}|>1\right),
\quad p_n =\textstyle\prod_{k=0}^{n-1} a_k,\notag\\
&\qquad d_j\,|\,\,p_n,\quad1\le j\le\sigma(p_n).\notag\\
b_1 &=1,\;
b_n = \textstyle\left\lceil\frac12\sum_{k=1}^{n-1}b_k\right\rceil.\notag\\
\label{p-t}\text{Then\ \ }a_n &=2^{b_n}, \quad \text{for\ }n>2.
\end{align}
\end{thm}

The first terms of $a_n$ and $b_n$ are \cite{S}\cite{Z}
\begin{align*}
a_{n\ge0} &=\{4,3,4,2,4,8,16,64,\ldots\},\quad
b_{n\ge1} =\{1,1,1,2,3,4,6,9,14\ldots\}.
\end{align*}
We need two lemmata.

\begin{lem}For $k>0$,
\begin{equation}\label{p-l1}
\sigma(3\cdot2^k)=2k+2.
\end{equation}
\end{lem}
\begin{proof}
This is true for $k=1$, and the set of divisors of~$3\cdot2^{k+1}$ is
the set of divisors of~$3\cdot2^k$ plus~$2^{k+1}$ and~$3\cdot2^{k+1}$
itself.
\end{proof}

\begin{lem}\label{l2}
Let $\delta(m)$ denote the smallest absolute value of the differences 
between complementary divisors of~$m>1$:
\begin{align}
\delta(m)\quad &= 
\quad\min\left(\,\left|\,d_j-\frac{m}{d_j}\,\right|\,\right),
\quad d_j|m,\quad1\le j\le\sigma(m).\notag\\
\intertext{Then}
\delta(3\cdot2^k)\quad &= \quad 2^{\lceil k/2\rceil},\quad k>0.\label{p-l2}
\end{align}
\end{lem}

\begin{proof}
Let us sort the divisors of $3\cdot2^k$ by size and call these~$D_j$:
$$D_{1\le j\le2k+2}=\{1,(2,3),\ldots,(2^i,\fract322^i),
(2^{i+1},\fract322^{i+1}),
\ldots,(2^k,\fract322^{k}),3\cdot2^k\}.
$$
Any smallest complementary divisor difference must be the one where
the divisors are in the exact middle of the sorted list, which, using
(\ref{p-l1}), is~$k+1$. And so, 
$\delta(3\cdot2^k)=D_{k+2}-D_{k+1}$.

Now, the proposition (\ref{p-l2}) is true for $k=1$. For every increase
of~$k$ by one, $\sigma(m)$ increases by two, and the index of the wanted 
pair of divisors 
increases by one, so $D_{k+2}-D_{k+1}$ goes through the values
\begin{align*}
\fract322^i-2^i &= 2^{i-1}\\
2^{i+1}-\fract322^i &= 2^{i-1}\\
\fract322^{i+1}-2^{i+1} &= 2^{i}\\
2^{i+2}-\fract322^{i+1} &= 2^{i}
\end{align*}
so it doubles every second step which is just the meaning of (\ref{p-l2}).
\end{proof}

Fixing the induction base at $\delta(p_3=48)=2^1=2^{b_3}$ to make sure
that $D_{k+2}-D_{k+1}>1$, the main proposition (\ref{p-t}) is now obvious, 
since the
powers of two in~$a_n$ behave the same way under multiplication as 
unity does in~$b_n$ under addition. 

\bigskip
Because the asymptotics of~$b_n$ are known\cite{C}, with
\begin{equation*}
b_n=\left\lceil c\cdot\left(\fract32\right)^n-\fract12\right\rceil,
\quad c=0.36050455619661495910154466\ldots,
\end{equation*}
the investigation of $a_{n\ge3}=2^{b_n}$ is settled, except for the
closed form for~$c$. Reble already proved\cite{R} that $b_n$ is
connected to the Josephus problem. Independently, our
numerics show that 
\begin{equation}
c=\frac29K(3),
\end{equation}
with $K(3)$ the universal constant in the same problem with~$q=3$,
a constant already discussed (\cite{OW}\cite{HH}),
and whose closed form is still unknown.

\end{document}